\documentclass[12pt]{amsart}

\usepackage{amsmath}

\theoremstyle{plain}
\newtheorem{theorem}{Theorem}
\newtheorem{corollary}{Corollary}
\newtheorem{conjecture}{Conjecture}
\newtheorem{lemma}{Lemma}

\theoremstyle{definition}

\theoremstyle{remark}

\numberwithin{equation}{section}
\numberwithin{table}{section}

\begin{document}
\title[Moments of the derivative]{Moments of the derivative of characteristic polynomials with an application to the Riemann zeta function}
\author {J.B. Conrey}
\address{American Institute of Mathematics,
360 Portage Ave, Palo Alto, CA 94306} \address{School of
Mathematics, University of Bristol, Bristol, BS8 1TW, United
Kingdom} \email{conrey@aimath.org}

\author{M.O. Rubinstein}
\address{Pure Mathematics, University of Waterloo, 
Waterloo, Ontario, N2L 3G1, Canada}
\email{mrubinst@uwaterloo.ca}

\author{N.C. Snaith}
\address{School of Mathematics,
University of Bristol, Bristol, BS8 1TW, United Kingdom}
\email{N.C.Snaith@bris.ac.uk}

\abstract
We investigate the moments of the derivative, on the unit circle, of characteristic polynomials of 
random unitary matrices and use this to formulate a conjecture for the moments of the derivative of 
the Riemann $\zeta$ function on the critical line. We do the same for the analogue of Hardy's 
$Z$-function, the characteristic polynomial multiplied by a suitable factor to make it real on the 
unit circle. Our formulae are expressed in terms of a determinant of a matrix whose entries involve 
the I-Bessel function and, alternately, by a combinatorial sum.
\endabstract

\maketitle

\section{Introduction}

Characteristic polynomials of unitary matrices serve as extremely
useful models for the Riemann zeta-function $\zeta(s)$. The
distribution of their eigenvalues give insight into the distribution
of zeros of the Riemann zeta-function and the values of these
characteristic polynomials give a model for the value distribution
of $\zeta(s)$. See the works [KS] and [CFKRS] for detailed
descriptions of how these models work. The important fact is that
formulas for the moments of the Riemann zeta-function are suggested
by the moments of the characteristic polynomials of unitary
matrices.

We consider two problems here: the moments of
the derivative of the characteristic polynomial $\Lambda_A(s)$ of an
$N \times N$ unitary matrix $A$, and
also the moments of the analogue of the derivative
of Hardy's $Z$-function, the characteristic polynomial multiplied by a suitable factor to
make it real on the unit circle.

In its simplest form our problem is to
give an exact formula, valid for complex $r$ with $\Re r>0$, of
the moments of the absolute value of the derivative of characteristic
polynomials
\begin{eqnarray}
    \int_{U(N)} |\Lambda_A'(1)|^r dA_N
\end{eqnarray}
or of
\begin{eqnarray}
    \int_{U(N)} |\mathcal Z_A'(1)|^r dA_N.
\end{eqnarray}
Here we are integrating against Haar measure on the unitary group, and
$\mathcal Z_A(s)$ is equal to $\Lambda_A(s)$ times a rotation factor that
makes it real on the unit circle. See the next section for the precise definition.

Unfortunately, we cannot yet solve either of these problems. However, we
can give asymptotic formulas when $r=2k$  for positive
integer values of $k$. The first two of these formulas involve the Maclaurin series
coefficients of a certain $k\times k$ determinant, while the third involves a
combinatorial sum.

\begin{theorem}

For fixed $k$ and $N\to \infty$ we have
\begin{eqnarray}
    \int_{U(N)} |\Lambda_A'(1)|^{2k} dA_N \sim b_k N^{k^2+2k},
\end{eqnarray}
where
\begin{eqnarray}
b_k = (-1)^{k(k+1)/2}
 \sum_{h=0}^k  \bigg( {k \atop h}\bigg)
    \bigg(\frac  d {dx}\bigg)^{k+h} \bigg(e^{-x} x^{-k^2/2}\det_{k\times k}
\big( I_{i+j-1}(2\sqrt{x}) \bigg)\bigg|_{ x=0},
\end{eqnarray}
and $I_\nu(z)$ denotes the modified Bessel function of the first kind.
\end{theorem}

\begin{theorem}

For fixed $k$ and $N\to \infty$ we have
\begin{eqnarray}
     \int_{U(N)} |\mathcal Z_A'(1)|^{2k} dA_N \sim b_k' N^{k^2+2k},
\end{eqnarray}
where
\begin{eqnarray}
b_k' = (-1)^{k(k+1)/2}   \bigg(\frac{d}{dx}\bigg)^{2k}\bigg(
 e^{-\tfrac x 2} x^{-k^2/2} \det_{k\times k}
\big(I_{i+j-1}(2\sqrt{x})\big)
\bigg)\bigg|_{x=0}.
\end{eqnarray}
\end{theorem}

We also have combinatorial description of $b_k'$.

\begin{theorem}
\begin{eqnarray}
b_k'=& \\
(-1)^{k(k+1)/2} &\sum_{m\in
P_O^{k+1}(2k)} \binom{2k}{m}\left(\frac{-1}2\right)^{m_0}
 \left( \prod_{i=1}^k \frac{1}{(2k-i+m_i)!} \right) \left(
\prod_{1\le i<j\le k} (m_j-m_i+i-j)\right), \notag
\end{eqnarray}
where $P_O^{k+1}(2k)$ denotes the set of partitions
$m=(m_0,\dots,m_k)$ of $2k$ into $k+1$ non-negative parts.

\end{theorem}

We have computed some values of $b_k$ and $b_k'$; these are tabulated at
the end of the paper.

Applying the random matrix theory philosophy suggests the conjecture:
\begin{conjecture}
\begin{equation}
    \frac{1}{T} \int_0^T |\zeta'(1/2+it)|^{2k} dt
    \sim a_k b_k \log(T)^{k^2+2k}
\end{equation}
and, similarly for Hardy's $Z$ function,
\begin{equation}
    \frac{1}{T} \int_0^T |Z'(1/2+it)|^{2k} dt
    \sim a_k b_k' \log(T)^{k^2+2k},
\end{equation}
where $a_k$ is the arithmetic factor
\begin{equation}
a_k =  \prod_p \left(1-\tfrac{1}{p}\right)^{k^2}
\sum_{m=0}^{\infty} \left( \frac{\Gamma(m+k)} {m! \Gamma(k)}
\right)^2 p^{-m}.
\end{equation}
\end{conjecture}

Remarks:
\begin{enumerate}

\item 
The factor $a_k$ is the same arithmetic contribution that arises
in the moments of the Riemann zeta function itself, see [KS] or
[CFKRS].  For an explanation of why these moments factor,
asymptotically, into the product of a contribution from the
primes, $a_k$, and a coefficient calculated via random matrix
theory, see [GHK]. 

\item 
In this paper we are only concerned with the leading asymptotics
of the moments of $\zeta'(1/2+it)$ and $Z'(1/2+it)$. Consequently, we use
the $k$-fold integrals for moments given below in Lemma~\ref{lemma:3}. If one
wishes to study lower order terms one would need to use the full moment conjecture
for $\zeta$ and $Z$ given in [CFKRS] as a $2k$ fold integral. 

\item Forrester and Witte have taken our Theorems 1-2 and managed to
find an alternate expression for $b_k$ and $b_k'$ involving a
Painlev\'e ${\rm III^{\prime}}\;$ equation,
and also an expression involving a certain generalised hypergeometric
function [FW, section 5].

\item
In his PhD thesis,
Chris Hughes gives a similar conjecture
for a more general problem involving mixed moments [Hug, conj 6.1].
What is new in our paper are the formulas for $b_k$ and $b_k'$.
For comparison, we state Hughes' formulation of the conjecture for
Hardy's $Z$ function. Let
\begin{eqnarray*}
I(h,k)=\int_0^T |Z(t)|^{2k-2h}\left|Z'(t)\right|^{2h}~dt.
\end{eqnarray*}
Hughes conjectures that
\begin{eqnarray*}
   I(h,k)\sim B(h,k) a_k f_k T(\log T)^{k^2+2h}
\end{eqnarray*}
where $a_k$ is given above, 
\begin{eqnarray*}
    f_k = \prod_{j=0}^{k-1}\frac{j!}{(k+j)!},
\end{eqnarray*}
and $B(h,k)$ is the constant that Hughes obtains
from the analogous moments for $\mathcal{Z}_A$
via random matrix theory:
\begin{eqnarray*}
   &&
   B(h,k)=\lim_{\beta\to 0}\frac{1}{\beta^{2h}}
   \sum_{n=0}^{2h} (-1)^{n-h}\left(
   {2h \atop n}\right) e^{-n\beta/2}
   \det\{b_{i,j}\},
\end{eqnarray*}
where
\begin{eqnarray*}
   b_{i,j}=\sum_{m=0}^{2h}
   \frac{(2k-n+i-1)!}{(2k-n+i-1+m)!}\left(
   {i+k-n-1+m \atop m}\right)
   \left({i+m-1\atop j-1}\right) \beta^m.
\end{eqnarray*}
To get Hughes' conjecture in this form, see (6.60), (6.51), and (6.52)
of [Hug] and replace $\beta$ by $\beta/(iN)$.

\item Brezin and Hikami [BH] attempt to use a similar approach to obtain a
theorem for moments of derivatives of characteristic polynomials, but there is an
error in their paper.

\item Numerically, we have observed that $b_k \sim 4^{-k} f_k$, as $k \to \infty$.
Other than
the power of 4, the r.h.s is the constant that appears in
the moments of characteristic polynomials [KS]:
\begin{eqnarray*}
    \int_{U(N)} |\Lambda_A(1)|^{2k} dA_N \sim f_k N^{k^2}, \quad \text{as $N \to \infty$}.
\end{eqnarray*}
A heuristic explanation is as follows. The large values of $|\Lambda_A'(1)|$ occur
near the large values of $|\Lambda_A(1)|$, namely when all the eigenvalues are close to $-1$.
The derivative of $\Lambda_A$ involves a sum of $N$ terms, each of 
which is missing, when the eigenvalues are close to $-1$, one factor of size roughly 2. 
A comparison with the $2k$-th moment
of $|\Lambda_A|$ thus gives an extra $N^{2k}$ and a factor of $2^{-2k}$.
We have not attempted to make this arguement rigorous.
We have also not  attempted to show that $b_k$ and $b_k'$ are non-zero.

\item The problem of moments of the derivative, are related, through
Jensen's formula, to the problem of zeros of the derivative. This
approach requires knowledge of the complex moments of the
derivative and we are only able to obtain integer moments.
For characteristic polynomials one is interested in studying the radial 
distribution of the zeros of the derivative. Francesco Mezzadri has the best 
results in this direction [Mez].
On the number theory side, one is interested in the horizontal distribution of
the zeros of $\zeta'$. Partial results have been obtained by
Levinson and Montgomery [LM], Conrey and Ghosh [CG], Soundararajan [Sou],
and Zhang [Z].

\end{enumerate}

\section{Notation}

If $A$ is an $N\times N$ matrix with complex entries $A=(a_{jk})$, we let $A^*$
be its conjugate transpose, i.e. $A^*=(b_{jk})$ where
$b_{jk}=\overline{a_{kj}}.$ $A$ is said to be unitary if $AA^*=I$.
We let $U(N)$ denote the group of all $N\times N$ unitary
matrices. This is a compact Lie group and has a Haar measure.

All of the eigenvalues of $A\in U(N)$ have absolute value 1; we write them
as
 \begin{eqnarray}e^{i\theta_1}, e^{i\theta_2}, \dots
e^{i\theta_N}.
\end{eqnarray}
The eigenvalues of $A^*$ are $e^{-i\theta_1},\dots, e^{-i\theta_N}$.
Clearly, the determinant, $\det A=\prod_{n=1}^N e^{i\theta_n}$ of a
unitary matrix is a complex number with absolute value equal to 1.

We are interested in computing various statistics about these eigenvalues.
Consequently, we identify all matrices in $U(N)$ that have the same
set of eigenvalues. The collection of matrices with the same set of
eigenvalues constitutes a conjugacy class in $U(N)$.
Weyl's integration formula [Weyl, pg 197] gives a simple way to perform 
averages over
$U(N)$ for functions $f$ that are constant on conjugacy classes.
Weyl's formula asserts that for such an $f$,
\begin{eqnarray}
    \int_{U(N)} f(A) ~d\mbox{Haar}=\int_{[0,2\pi]^N}
    f(\theta_1,\dots,\theta_N)dA_N
\end{eqnarray}
where
\begin{eqnarray}
    dA_N =\prod_{1\le j<k\le
    N}\big|e^{i\theta_k}-e^{i\theta_j}\big|^2 ~\frac{d\theta_1 \dots
    d\theta_N}{N! (2\pi)^N}.
\end{eqnarray}

The characteristic polynomial of a matrix $A$ is denoted
$\Lambda_A(s)$ and is defined by
\begin{eqnarray}\Lambda_A(s)=\det(I-sA^*)=\prod_{n=1}^N(1-se^{-i\theta_n}).\end{eqnarray}
The roots of $\Lambda_A(s)$ are the eigenvalues of $A$ and are on the unit
circle. Notice that this definition of the characteristic polynomial
differs slightly from the usual definition in that it has an extra factor of
$\det(A^*)$.
We regard $\Lambda_A(s)$ as an analogue of the Riemann zeta-function
and this definition is chosen so as to resemble the
Hadamard product of $\zeta$.

The characteristic polynomial  satisfies the functional equation
\begin{eqnarray} \Lambda_A(s)&=&(-s)^N\prod_{n=1}^N
e^{-i\theta_n}\prod_{n=1}^N
(1-e^{i\theta_n}/s)\\
&=&(-1)^N \det A^* ~s^N~\Lambda_{A^*}(1/s).\end{eqnarray}

We define the $\mathcal Z$-function by
\begin{eqnarray}\mathcal Z_A(s)=e^{-\pi i N/2}  e^{i\sum_{n=1}^N
\theta_n/2}
s^{-N/2}\Lambda_A(s);
\end{eqnarray}
here if $N$ is odd,
 we use the branch of the square-root function that is positive
for positive real $s$.
The functional equation for $\mathcal Z$ is
\begin{eqnarray} \mathcal Z_A(s)=(-1)^N \mathcal Z_{A^*}(1/s) .
\end{eqnarray}
Note that
\begin{eqnarray}
\overline{\mathcal Z_A(e^{i\theta})}=\mathcal Z_{A }(e^{i\theta})
\end{eqnarray}
so that $\mathcal Z_A (e^{i\theta})$ is real when $\theta$ is real.
We regard $\mathcal Z_A(e^{i\theta})$ as an analogue of Hardy's function
$Z(t)$.

We let $I_n$ be the usual modified Bessel function with power series
expansion
\begin{eqnarray}
\label{eq:besselpower}
 I_n(x)=\big(\frac x 2
\big)^n\sum_{j=0}^\infty \frac {x^{2j}}{2^{2j}(n+j)! j!}.
\end{eqnarray}

The way that the I-Bessel function enters our calculation is through
the following formula:
\begin{eqnarray}
\frac{1}{2\pi i}\int_{|z|=1} \frac {e^{Lz+t/z}}{z^{2k}} ~dz =
\frac{L^{2k-1}
I_{2k-1}(2\sqrt{Lt})}{(Lt)^{k-1/2}}.
\end{eqnarray}
This formula can be proven by
comparing the coefficient of $z^{2k-1}$ in $e^{Lz+t/z}$ with the power
series formula for $I_{2k-1}$.

We let $\Delta(z_1,\dots,z_k)$ denote the Vandermonde determinant
\begin{eqnarray}
\Delta(z_1,\dots,z_k)=\det_{k\times k}\big( z_i^{j-1}\big).
\end{eqnarray}
We often omit the subscripts and write $\Delta(z)$ in place of
$\Delta(z_1,\dots,z_k)$.
Also, we allow differential operators as the arguments, such as
\begin{eqnarray}
    \Delta\bigg(\frac{d}{dL}\bigg)
    =\Delta\bigg(\frac{d}{dL_1},\dots,\frac{d}{dL_k} \bigg)
    = \det_{k\times k}\bigg( \bigg(\frac{d}{dL_i}\bigg)^{j-1}\bigg).
\end{eqnarray}
The key fact about the Vandermonde is that
\begin{eqnarray} \Delta(z_1,\dots,z_k)=\prod_{1\le i<j\le k}(z_j-z_i).
\end{eqnarray}

We let
\begin{eqnarray}
z(x)=\frac{1}{1-e^{-x}}=\frac{1}{x}+O(1).
\end{eqnarray}
The function $z(x)$ plays the role for random matrix theory that
$\zeta(1+x)$ plays in the theory of moments of the Riemann
zeta-function. See for example pages 371--372 of [CFKRS2].

We let $\Xi$ denote the subset of permutations $\sigma\in S_{2k}$ of
$\{1,2,\dots,2k\}$ for which
\begin{eqnarray}
\sigma(1)<\sigma(2)<\dots < \sigma(k)\end{eqnarray}
and
\begin{eqnarray}
\sigma(k+1)<\sigma(k+2)<\dots < \sigma(2k).\end{eqnarray}

We let $P_O^{k+1}(2k)$ be the set of partitions
$m=(m_0,\dots,m_k)$ of $2k$ into $k+1$ non-negative parts. This quantity
arises
from the multinomial expansion
\begin{eqnarray}
\label{eq:multinom}
 (x_0+x_1+\cdots+x_k)^{2k}=\sum_{m\in P_O^{k+1}(2k)}
\binom{2k}{m} x_0^{m_0}\cdots x_k^{m_k}
\end{eqnarray}
where
\begin{eqnarray}
\binom{2k}{m}=\frac{(2k)!}{m_0!\dots m_k!}.
\end{eqnarray}

\section{Lemmas}

The main tool in proving theorems 1-2 is to take formulas (Lemma
3) for moments of characteristic polynomials with shifts,
differentiate these with respect to the shifts, and then set the
shifts equal to zero. This gives $k$-fold contour integrals. To
separate the integrals involved, we introduce extra parameters and
differential operators to pull out a portion of these integrands.

\begin{lemma}
Assume that $\alpha_1,\ldots,\alpha_{2k}$ are distinct complex numbers.
We have
\begin{eqnarray}
\label{eq:lem1}
    && \int_{U(N)} \prod_{j=1}^k \Lambda_{A}(e^{-\alpha_j})
    \Lambda_{A^*}
    (e^{\alpha_{j+k}})
    ~dA_N =
    \sum_{\sigma\in \Xi}
    e^{N\sum_{j=1}^k(\alpha_{\sigma(j)}-\alpha_{j } )}
    \prod_{1\le i,j\le k}z(\alpha_{\sigma ( j)}-\alpha_{\sigma(k+i)}).
\end{eqnarray}
\end{lemma}
This is proven in section 2 of [CFKRS2]. See formulas (2.5), (2.16), and (2.21)
of that paper. The definition given there of the
characteristic polynomial differs slightly from the one we use here, and that
introduces some extra exponential factors in (2.21) of the aforementioned paper,
and also necessitates replacing the $\alpha$'s by $-\alpha$'s.

Since
\begin{eqnarray}
\mathcal Z_A(e^{-\alpha_j})\mathcal Z_{A^*}(e^{\alpha_{j+k}})
=(-1)^N e^{N(\alpha_{j }-\alpha_{j+k})/2}\Lambda_A(e^{-\alpha_j})
\Lambda_{A^*}(e^{\alpha_{j+k}})
\end{eqnarray}
we can write a corresponding lemma for $\mathcal Z$.

\begin{lemma}
Assume that $\alpha_1,\ldots,\alpha_{2k}$ are distinct complex numbers. Then
\begin{eqnarray}
    && \int_{U(N)} \prod_{j=1}^k \mathcal Z_A(e^{-\alpha_j})
    \mathcal Z_{A^*}
    (e^{\alpha_{j+k}})
    ~dA_N\\
    &&\qquad
    =(-1)^{Nk} e^{-\tfrac N2\sum_{j=1}^{2k}\alpha_j}\sum_{\sigma\in \Xi}
    e^{N\sum_{j=1}^k\alpha_\sigma(j) }
    \prod_{1\le i,j\le k}z(\alpha_{\sigma ( j)}-\alpha_{\sigma(k+i)}). \notag
\end{eqnarray}
\end{lemma}

We can express the sums in the last two lemmas as integrals. Thus we have

\begin{lemma}
\label{lemma:3}
Assume that all of the $\alpha_j$ are smaller
than 1 in absolute value. Then
\begin{eqnarray}
\label{eq:lem3a}
    && \int_{U(N)}
    \prod_{j=1}^k \Lambda_{A}(e^{-\alpha_j})
    \Lambda_{A^*}
    (e^{\alpha_{j+k}})
    ~dA_N
    \\
    &&\qquad
    =
    \frac{1 }{k!(2\pi i)^k}\int_{|w_i|=1}
    e^{N\sum_{j=1}^k (w_j-\alpha_j)}
    \prod_{1\le i\le k\atop 1\le j\le 2k}
    z(w_i-\alpha_j)\prod_{i\ne j}z(w_i-w_j)^{-1}\prod_{j=1}^k dw_j \notag
\end{eqnarray}
and
\begin{eqnarray}
    && \int_{U(N)} \prod_{j=1}^k \mathcal Z_A(e^{-\alpha_j})
    \mathcal Z_{A^*}
    (e^{\alpha_{j+k}})
     ~dA_N\\
     &&\qquad
    = (-1)^{Nk} \frac{e^{-\tfrac N2\sum_{j=1}^{2k}\alpha_j}}{k!(2\pi i)^k}\int_{|w_i|=1}
    e^{N\sum_{j=1}^k w_j}
    \prod_{1\le i\le k\atop 1\le j\le 2k}
    z(w_i-\alpha_j)\prod_{i\ne j}z(w_i-w_j)^{-1}\prod_{j=1}^k dw_j. \notag
\end{eqnarray}
\end{lemma}
In this lemma, and its corollary below, we do not require the
$\alpha_j's$ to be distinct.   The proof of Lemma \ref{lemma:3} is
a straight-forward evaluation of the residues in the integral in
(\ref{eq:lem3a}), arising from the factor $z(w_i-\alpha_j)$, to
obtain the $\binom{2k}{k}$ terms in (\ref{eq:lem1}).  Each of the
$k$ integrals in (\ref{eq:lem3a}) results in a sum over $2k$
residues, but due to the factor $\prod_{i\neq j} z(w_i-w_j)^{-1}$,
any one of these $2k^2$ terms is zero if the residue of two of the
integrals, say $w_i$ and $w_j$, are evaluated at the same point
$\alpha_{\ell}$.

Using the fact that $z(w)=1/w+O(1)$ we easily deduce
\begin{corollary}
Suppose that $\alpha_j=\alpha_j(N)$ and $|\alpha_j| \ll 1/N$ for each $j$. Then
\begin{eqnarray}
    && \int_{U(N)}
    \prod_{j=1}^k \Lambda_{A }(e^{-\alpha_j})
    \Lambda_{A^*}
    (e^{\alpha_{j+k}})
    ~dA_N
    \\
    &&\qquad
    =
    \frac{1}{k!(2\pi i)^k}\int_{|w_i|=1}
    e^{N\sum_{j=1}^k(w_j-\alpha_j)}
    \frac{\prod_{i\ne j}(w_i-w_j)}
    {\prod_{1\le i\le k\atop 1\le j\le 2k}
    (w_i-\alpha_j) }\prod_{j=1}^k dw_j +O(N^{k^2-1}) \notag
\end{eqnarray}
with an implicit  constant independent of $N$; similarly,
\begin{eqnarray}
\label{eq:coroll2}
   && \int_{U(N)} \prod_{j=1}^k \mathcal Z_A(e^{-\alpha_j})
   \mathcal Z_{A^*}
   (e^{\alpha_{j+k}})
   ~dA_N \\
   &&\qquad
   = (-1)^{Nk} \frac{e^{-\tfrac N2\sum_{j=1}^{2k}\alpha_j}}{k!(2\pi i)^k}\int_{|w_i|=1}
   e^{N\sum_{j=1}^k  w_j }
   \frac{\prod_{i\ne j}(w_i-w_j)}
   {\prod_{1\le i\le k\atop 1\le j\le 2k}
   (w_i-\alpha_j) }\prod_{j=1}^k dw_j +O(N^{k^2-1}) \notag
\end{eqnarray}
\end{corollary}

\begin{lemma}
Let $f$ be $k-1$ times differentiable, $k \geq 1$. Then
\label{lemma:single vand}
$$
    \Delta(\frac{d}{dL})
    \prod_{i=1}^k f(L_i)
    = \det_{k\times k} \left( f^{(j-1)}(L_i) \right)
$$
where by $\Delta(d/dL)$ we mean the differential operator
\begin{eqnarray}
    \prod_{1\le i<j\le k}\bigg(\frac
    {d}{dL_j}-\frac{d}{dL_i}\bigg)
    =\det_{k\times k} \left( \frac {d^{j-1}}{dL_i^{j-1}}\right).
\end{eqnarray}
\end{lemma}

\vspace{.2cm}
\noindent{\bf Proof.}
This follows using the definition of the Vandermonde determinant.
Noticing that row $i$ of the matrix only involves $L_i$, we factor
the product into the determinant.

\begin{lemma}  Let $f$ be $2k-2$ times differentiable. Then
\label{lemma:two vandermondes}
\begin{eqnarray}
    \Delta^2\bigg(\frac d{dL}\bigg)\big(\prod_{i=1}^k
    f(L_i)\big)\bigg|_{L_i=L}=k! \det_{k\times k}
     \big(f^{(i+j-2)}(L)\big).
\end{eqnarray}
More generally, suppose that
$g(L_1,\dots,L_k)=\sum_{r=1}^Ra_r \prod_{i=1}^k f_{r,i}(L_i)$ is a
symmetric function
of its $k$ variables. Then
\begin{eqnarray}
\Delta^2\bigg(\frac d{dL}\bigg)
g(L_1,\dots,L_K)\bigg|_{L_j=L}=k! \sum_{r=1}^Ra_r \det_{k\times k}
 \big(f_{r,i}^{(i+j-2)}(L)\big).
\end{eqnarray}
\end{lemma}

\vspace{.2cm}
\noindent{\bf Proof.}
Applying the Vandermonde a second time to Lemma~\ref{lemma:single vand} we get
\begin{equation}
    \Delta\left(\frac{d}{dL}\right)
    \det_{k\times k} \left( f^{(j-1)}(L_i) \right).
\end{equation}
Expand the determinant as a sum over all permutations $\mu$ of the numbers
$1,2,\ldots,k$:
\begin{equation}
    \det_{k\times k} \left( f^{(j-1)}(L_i) \right)
    = \sum_{\mu} \text{sgn}(\mu) \prod_{i=1}^k f^{\mu_i-1} (L_i).
\end{equation}
Apply Lemma~\ref{lemma:single vand} to find that a typical term above equals
\begin{equation}
    \text{sgn}(\mu)
    \det_{k\times k}
    \big(f^{(\mu_i+j-2)}(L_i)\big).
\end{equation}
Setting $L_i=L$ for $1 \leq i \leq k$, we may rearrange the rows so as to
undo the permutation $\mu$. This introduces another $\text{sgn}(\mu)$ in front of the determinant
and gives
\begin{equation}
    \det_{k\times k}
    \big(f^{(i+j-2)}(L)\big).
\end{equation}
Since there are $k!$ permutations $\mu$, we get
\begin{equation}
    k! \det_{k\times k}
    \big(f^{(i+j-2)}(L)\big).
\end{equation}
The proof of the second part of the lemma is left to the reader.

\begin{lemma}
Suppose that $P$ and $Q$ are polynomials with
$Q(w)=\prod_{j=1}^{2k}(w-\alpha_j)$
and $\max |\alpha_j| <c.$
Then
\begin{eqnarray}
    \frac{1}{2\pi i} \int_{|w|=c} \frac{e^{wL}}{w}\frac{P(w)}{Q(w)}
    dw
     =P\bigg(\frac{d}{dL}\bigg) \int_{\sum_{j=1}^{2k}x_j \leq L}
    e^{\sum_{j=1}^{2k} x_j \alpha_j} \prod_{j=1}^{2k}dx_j.
\end{eqnarray}
\end{lemma}

\vspace{.2cm} \noindent{\bf Proof.} Since
\begin{equation}
P\bigg(\frac{d}{dL}\bigg)e^{wL}=e^{wL}P(w),
\end{equation}
the derivatives can be pulled outside the integral immediately.
With the Laplace transform pair $e^{x\alpha}$ and
$\frac{1}{w-\alpha}$, related by
\begin{equation}
e^{x\alpha}=\frac{1}{2\pi i} \int_{|w|=c} \frac{e^{wx}}{w-\alpha}
dw,
\end{equation}
we merely apply repeatedly the Laplace convolution formula, which
for Laplace transform pairs $f_i$ and $\phi_i$ states that
\begin{equation}
\frac{1}{2\pi i} \int_{|w|=c}
\phi_1(s)\phi_2(s)e^{sx}ds=\int_0^{x} f_1(y)f_2(x-y) dy,
\end{equation}
to evaluate the Laplace transform of the product
$\frac{1}{w\;\prod_{j=1}^{2k}(w-\alpha_j)}$.

\begin{lemma}
We have
\begin{eqnarray}
\int_{\sum_{j=1}^{2k}x_j\le L} x_1\dots x_n \prod_{j=1}^{2k}dx_j
=\frac{L^{2k+n}}{(2k+n)!}.
\end{eqnarray}
\end{lemma}

This lemma can be proved in a straight-forward manner by
induction.

\section{Proofs}

We now give the proofs of our identities for the leading terms of the
moments
of the derivatives of $\Lambda$ and $\mathcal Z$. We begin with the proof
of Theorem 2 for
$\mathcal Z$ as it is slightly easier.

\vspace{.2cm}
\noindent{\bf Proof of Theorem 2.}
A differentiated form of the second formula of Corollary 1 gives us
\begin{eqnarray}
    \label{eq:test}
    &&\prod_{j=1}^{2k} \frac{d}{d\alpha_j} \int_{U(N)}\prod_{h=1}^k
    \mathcal Z_A(e^{-\alpha_h})\mathcal Z_{A^*}
    (e^{\alpha_{k+h}})dA_N
    = \\
    &&(-1)^{\frac{k(k-1)}{2}+kN}
    \prod_{j=1}^{2k} \frac{d}{d\alpha_j}\frac{e^{-\tfrac
    N2\sum_{j=1}^{2k}\alpha_j}}
    {k!(2\pi i)^k}\int_{|w_i|=1}
    e^{N\sum_{j=1}^k w_j }
    \frac{ \Delta^2(w)}
    {\prod_{1\le i\le k\atop 1\le j\le 2k}
    (w_i-\alpha_j) }\prod_{j=1}^k dw_j \notag \\
    &&+O(N^{k^2+2k-1}), \notag
\end{eqnarray}
provided that $\alpha_j=\alpha_j(N)\ll 1/N$. Notice that
\begin{eqnarray}
\label{eq:derivZA}
   \frac{d}{d\alpha} \mathcal{Z}_A(e^{-\alpha})\big|_{\alpha=0} =
   -\frac{d}{ds} \mathcal{Z}_A(s)\big|_{s=1}=-\mathcal{Z}'_A(1)
\end{eqnarray}
and
\begin{eqnarray}
\label{eq:derivZA*}
     \frac{d}{d\alpha}
     \mathcal{Z}_{A^*}(e^{\alpha})\big|_{\alpha=0} =
     \mathcal{Z}'_{A^*}(1)=(-1)^N\overline{\mathcal{Z}'_A (1)}.
\end{eqnarray}
So,
\begin{eqnarray}
    \label{eq: Z' as contour}
    &&\int_{U(N)} |\mathcal Z_A'(1)|^{2k} dA_N \\
    &&=
    (-1)^{\frac{k(k+1)}{2}}
    \prod_{j=1}^{2k} \frac{d}{d\alpha_j}\frac{e^{-\tfrac
    N2\sum_{j=1}^{2k}\alpha_j}}
    {k!(2\pi i)^k}\int_{|w_i|=1}
    e^{N\sum_{j=1}^k w_j }
    \frac{ \Delta^2(w)}
    {\prod_{1\le i\le k\atop 1\le j\le 2k}
    (w_i-\alpha_j) }\prod_{j=1}^k dw_j \bigg|_{\alpha=0} \notag \\
    &&+O(N^{k^2+2k-1}). \notag
\end{eqnarray}
The sign here arises as the $(-1)^{kN}$ from (\ref{eq:derivZA*})
cancels the same factor in (\ref{eq:test}), we have a $(-1)^k$
from (\ref{eq:derivZA}) and we pick up the
$(-1)^{\frac{k(k-1)}{2}}$ in (\ref{eq:test}) through writing the
factor $\prod_{i\neq j} (w_i-w_j)$ in (\ref{eq:coroll2}) as
$\Delta^2(w)$ above.

 To separate the integrals, we introduce extra
parameters $L_i$ and move the Vandermonde polynomial outside the
integral as a differential operator, getting
\begin{eqnarray}
    \label{eq:Z' vandermonde outside}
    &&(-1)^{\frac{k(k+1)}{2}}
    \prod_{j=1}^{2k} \frac{d}{d \alpha_j}
    \frac{ \Delta^2(d/dL) e^{-\tfrac N2\sum_{j=1}^{2k}\alpha_j} }
    {k!(2\pi i)^k}
    \int_{|w_j|=1}
    \frac{e^{\sum_{i=1}^k L_iw_i} }{\prod_{1\le i\le k\atop 1\le j\le 2k}
    (w_i-\alpha_j) }\prod_{j=1}^k dw_j  \bigg|_{\alpha=0,L_i=N} \\
    &&+O(N^{k^2+2k-1}). \notag
\end{eqnarray}

Next, we observe that
\begin{eqnarray}
    \frac{d}{d \alpha } \frac{e^{-\tfrac N2\alpha}}{\prod_{1\le i\le k }
    (w_i-\alpha) }\bigg|_{\alpha=0}
    =\frac{1}{\prod_{i=1}^kw_i }\bigg(\sum_{j=1}^{k}\frac{1}{w_j}-\frac
    N2\bigg)
\end{eqnarray}
so that (\ref{eq:Z' vandermonde outside}) equals, without the $O$ term,
\begin{eqnarray}
    &&(-1)^{\frac{k(k+1)}{2}}
    \frac{  \Delta^2(d/dL)}{k!(2\pi i)^k}\int_{|w_j|=1}
    \frac{e^{\sum_{i=1}^k L_iw_i}\bigg(\sum_{j=1}^{k}\frac{1}{w_j}-\frac
    N2\bigg)^{2k} }
    {\prod_{i=1}^k
    w_i^{2k} }\prod_{j=1}^k dw_j \bigg|_{ L_i=N}.
\end{eqnarray}
Introducing another auxiliary variable $t$, this can be expressed as
\begin{eqnarray}
    (-1)^{\frac{k(k+1)}{2}}
    \frac{ \Delta^2(d/dL)\big(\frac{d}{dt}\big)^{2k}e^{-Nt/2}}
    {k!(2\pi i)^k}\int_{|w_j|=1}
    \frac{e^{\sum_{i=1}^k L_iw_i+t/w_i} }{\prod_{i=1}^k
    w_i^{2k} }\prod_{j=1}^k dw_j \bigg|_{ L_i=N,t=0}.
\end{eqnarray}
This allows us to separate the integrals and we get
\begin{eqnarray}
    (-1)^{\frac{k(k+1)}{2}}
    \frac{  \Delta^2(d/dL)\big(  d/dt\big)^{2k}e^{-Nt/2}}{k!}\prod_{i=1}^k
    \bigg(\frac{1}{2\pi i}\int_{|w|=1}\frac{e^{  L_iw +t/w  }}{
    w ^{2k} }dw\bigg)\bigg|_{L_i=N,t=0} .
\end{eqnarray}
The integral  evaluates to
\begin{eqnarray}
    \frac{L_i^{2k-1}
    I_{2k-1}(2\sqrt{L_it})}{(L_it)^{k-1/2}}
\end{eqnarray}
as noted earlier.  Thus,
\begin{eqnarray}
     \label{eq:Z' bessel}
     &&\int_{U(N)} |\mathcal Z_A'(1)|^{2k} dA_N  = \\
     && (-1)^{\frac{k(k+1)}{2}}
     \frac{ \Delta(d/dL)\big(  d/dt\big)^{2k}e^{-Nt/2}}
     {k!}\bigg(\prod_{i=1}^k
     \frac{L_i^{2k-1}
     I_{2k-1}(2\sqrt{L_it})}{(L_it)^{k-1/2}}\bigg)\bigg|_{L_i=N,t=0} \notag \\
     &&+O(N^{k^2+2k-1}). \notag
\end{eqnarray}
So, letting
\begin{eqnarray}
    f_t(L_i)=\frac{L_i^{2k-1}
    I_{2k-1}(2\sqrt{L_it})}{(L_it)^{k-1/2}},
\end{eqnarray}
we have, by Lemma~\ref{lemma:two vandermondes}, that (\ref{eq:Z' bessel}) equals
\begin{eqnarray}
    \label{eq:Z' bessel det 1}
    (-1)^{\frac{k(k+1)}{2}}
    \bigg( \frac{ d}{dt}\bigg)^{2k}e^{-Nt/2}
    \bigg(\det_{k\times k}
    \big(f_t^{(i+j-2)}(N)
    \bigg)\bigg|_{ t=0} +O(N^{k^2+2k-1}).
\end{eqnarray}

Now we see, from (\ref{eq:besselpower}), that
\begin{eqnarray}
    f_t(L)=\sum_{r=0}^\infty
    \frac{t^rL^{2k-1+r}}{r!(2k-1+r)!},
\end{eqnarray}
so that if $\mu\le 2k-1$, then
\begin{eqnarray}
    f^{(\mu)}(L)=\sum_{r=0}^\infty \frac{t^rL^{2k-1-\mu+r}}{r!(2k-1-\mu+r)!}
    =\bigg(\frac{L}{t}\bigg)^{(2k-1-\mu)/2}I_{2k-1-\mu}(2\sqrt{Lt}).
\end{eqnarray}
Therefore, (\ref{eq:Z' bessel det 1}) equals
\begin{eqnarray}
\label{eq:Ideta}
    &&(-1)^{\frac{k(k+1)}{2}}
    \bigg( \frac{ d}{dt}\bigg)^{2k} e^{-Nt/2}
    \det_{k\times k}
    \bigg( \bigg(\frac Nt\bigg)^{(2k+1-i-j)/2}I_{2k+1-i-j}(2\sqrt{Nt})
    \bigg)\bigg|_{ t=0} \notag \\
    &&+O(N^{k^2+2k-1}).
\end{eqnarray}
Clearly $\det_{k}(a_{i,j})=\det_{k}(a_{k+1-i,k+1-j})$, therefore
(\ref{eq:Ideta}) can be written as
\begin{eqnarray}
    (-1)^{\frac{k(k+1)}{2}}
    \bigg( \frac{ d}{dt}\bigg)^{2k} e^{-Nt/2}
    \det_{k\times k}
    \bigg( \bigg(\frac Nt\bigg)^{(i+j-1)/2}I_{i+j-1}(2\sqrt{Nt})
    \bigg)\bigg|_{ t=0} .
\end{eqnarray}
If we substitute $x=Nt$, then $d/dt=N d/dx$ and we get
\begin{eqnarray}
    &&(-1)^{\frac{k(k+1)}{2}}
    N^{2k} \bigg( \frac{ d}{dx}\bigg)^{2k} e^{-x/2}
    \det_{k\times k}
    \bigg( \bigg(\frac {N^2}x\bigg)^{(i+j-1)/2}I_{i+j-1}(2\sqrt{x})
    \bigg)\bigg|_{ x=0} \\
    =&&
    (-1)^{\frac{k(k+1)}{2}}
    N^{k^2+2k} \bigg(\frac  d {dx}\bigg)^{2k} \bigg(e^{-x/2}
    x^{-k^2/2}\det_{k\times k}
    \big(   I_{i+j-1}(2\sqrt{x})
    \bigg)\bigg|_{ x=0} ,\notag
\end{eqnarray}
since $\det_k (M^{i+j-1}a_{i,j})=M^{k^2}\det_k(a_{i,j})$
as is seen by factoring $M^j$ out of the $j$th column and $M^{i-1}$ out of
the $i$th row.
This completes the proof of Theorem 2.

\vspace{.2cm}
\noindent{\bf Proof of Theorem 1.}
Turning to Theorem 1's proof, we begin as before, but with
a differentiated form of the first formula of Corollary 1:
\begin{eqnarray}
    \label{eq:lambda contour}
    &&\prod_{j=1}^{2k} \frac{d}{d\alpha_j} \int_{U(N)}\prod_{h=1}^k
    \Lambda_A(e^{-\alpha_h})\Lambda_{A^*}
    (e^{\alpha_{k+h}})dA_N \\
    &&=
    (-1)^{\frac{k(k-1)}{2}}
    \prod_{j=1}^{2k} \frac{d}{d\alpha_j}\frac{1}{k!(2\pi i)^k}\int_{|w_i|=1}
    e^{N\sum_{j=1}^k(w_j-\alpha_j)}
    \frac{\Delta^2(w)}
    {\prod_{1\le i\le k\atop 1\le j\le 2k}
    (w_i-\alpha_j) }\prod_{j=1}^k dw_j \notag \\
    &&+O(N^{k^2+2k-1}), \notag
\end{eqnarray}
provided that $\alpha_j\ll 1/N$. Now
\begin{eqnarray}
\frac{d}{d\alpha} \Lambda_A(e^{-\alpha})\big|_{\alpha=0} =
-\frac{d}{ds} \Lambda_A(s)\big|_{s=1} = -\Lambda'_A(1)
\end{eqnarray}
and
\begin{eqnarray}
 \frac{d}{d\alpha} \Lambda_{A^*}(e^{\alpha})\big|_{\alpha=0} =    \Lambda'_{A^*}(1)=
     \overline{\Lambda'_A(1)},
\end{eqnarray}
hence setting $\alpha=0$, (\ref{eq:lambda contour}) becomes
\begin{eqnarray}
    \label{eq:lambda contour 2}
    && \int_{U(N)} |\Lambda_A'(1)|^{2k} dA_N \\
    &&=
    (-1)^{\frac{k(k+1)}{2}}
    \prod_{j=1}^{2k} \frac{d}{d\alpha_j}\frac{1}{k!(2\pi i)^k}\int_{|w_i|=1}
    e^{N\sum_{j=1}^k(w_j-\alpha_j)}
    \frac{\Delta^2(w)}
    {\prod_{1\le i\le k\atop 1\le j\le 2k}
    (w_i-\alpha_j) }\prod_{j=1}^k dw_j \bigg|_{\alpha_j=0} \notag \\
    &&+O(N^{k^2+2k-1}). \notag
\end{eqnarray}
Introducing variables $L_i$ as before, the above equals, without the $O$ term
\begin{eqnarray}
    (-1)^{\frac{k(k+1)}{2}}
    \frac{\prod_{j=1}^{2k}(d/d \alpha_j) \Delta^2(d/dL)}
    {k!(2\pi i)^k}\int_{|w_j|=1}
    \frac{e^{\sum_{i=1}^k (L_iw_i-N\alpha_i)} }{\prod_{1\le i\le k\atop 1\le
    j\le 2k}
    (w_i-\alpha_j) }\prod_{j=1}^k dw_j  \bigg|_{\alpha_j=0,L_i=N}.
\end{eqnarray}
Performing the differentiations with respect to the $\alpha_j$ leads us to
\begin{eqnarray}
    \label{eq:another formula}
    \\
    (-1)^{\frac{k(k+1)}{2}}
    \frac{  \Delta^2(d/dL)}{k!(2\pi i)^k}\int_{|w_j|=1}
    \frac{e^{\sum_{i=1}^k
    L_iw_i}\bigg(\sum_{j=1}^{k}\frac{1}{w_j}-N\bigg)^{k}
    \bigg(\sum_{j=1}^{k}\frac{1}{w_j} \bigg)^{k}}
    {\prod_{i=1}^k
    w_i^{2k} }\prod_{j=1}^k dw_j \bigg|_{ L_i=N}. \notag
\end{eqnarray}
Now we write
\begin{eqnarray}
    \bigg(\sum_{j=1}^{k}\frac{1}{w_j}-N\bigg)^{k}
    \bigg(\sum_{j=1}^{k}\frac{1}{w_j} \bigg)^{k}
    &=& \bigg(\sum_{j=1}^{k}\frac{1}{w_j}-N\bigg)^{k}
    \bigg(\sum_{j=1}^{k}\frac{1}{w_j}-N+N \bigg)^{k}\\
    &=& \sum_{h=0}^k \bigg( {k \atop h}\bigg)N^{k-h}
    \bigg(\sum_{j=1}^{k}\frac{1}{w_j}-N\bigg)^{k+h}. \notag
\end{eqnarray}

Introducing the auxiliary variable $t$, (\ref{eq:another formula}) can be expressed as
\begin{eqnarray}
    \\
    &&(-1)^{\frac{k(k+1)}{2}}
    \sum_{h=0}^k  \bigg( {k \atop h}\bigg)N^{k-h}
    \frac{  \Delta^2(d/dL)\big(\frac{d}{dt}\big)^{k+h}e^{-Nt}}
    {k!(2\pi i)^k}\int_{|w_j|=1}
    \frac{e^{\sum_{i=1}^k L_iw_i+t/w_i} }{\prod_{i=1}^k
    w_i^{2k} }\prod_{j=1}^k dw_j \bigg|_{ L_i=N,t=0} \notag \\
    &&=(-1)^{\frac{k(k+1)}{2}}\sum_{h=0}^k\binom{k}{h} N^{k-h} \frac{  \Delta^2(d/dL)\big(  d/dt\big)^{k+h}e^{-Nt}}{k!}\prod_{i=1}^k
    \bigg(\frac{1}{2\pi i}\int_{|w|=1}\frac{e^{  L_iw +t/w  }}{
    w ^{2k} }dw\bigg)\bigg|_{L_i=N,t=0} . \notag
\end{eqnarray}

Proceeding as before we arrive at
\begin{eqnarray}
    && \int_{U(N)} |\Lambda_A'(1)|^{2k} dA_N \\
    &&=
    (-1)^{\frac{k(k+1)}{2}}
    N^{k^2+2k}\sum_{h=0}^k  \bigg( {k \atop h}\bigg)
    \bigg(\frac  d {dx}\bigg)^{k+h} \bigg(e^{-x} x^{-k^2/2}\det_{k\times k}
    \big(   I_{i+j-1}(2\sqrt{x})
    \bigg)\bigg|_{ x=0} \notag \\
    &&+O(N^{k^2+2k-1}). \notag
\end{eqnarray}

\vspace{.2cm}
\noindent{\bf Proof of Theorem 3.}
We now give the proof of Theorem 3. We rewrite
equation~(\ref{eq:test}) as
\begin{eqnarray}
    &&\prod_{j=1}^{2k} \frac{d}{d\alpha_j} \int_{U(N)}\prod_{h=1}^k
    \mathcal Z_A(e^{-\alpha_h})\mathcal Z_{A^*}dA_N
    = \\
    &&(-1)^{\frac{k(k-1)}{2}+kN}
    \prod_{j=1}^{2k} \frac{d}{d\alpha_j}\frac{e^{-\tfrac
    N2\sum_{j=1}^{2k}\alpha_j}}
    {k!(2\pi i)^k}\int_{|w_i|=1}
    e^{N\sum_{i=1}^k w_i }
    \frac{ \Delta^2(w)\prod_{i=1}^k w_i}
    {\prod_{1\le i\le k\atop 1\le j\le 2k}
    (w_i-\alpha_j) }\prod_{i=1}^k \frac{dw_i}{w_i} \notag \\
    &&+O(N^{k^2+2k-1}). \notag
\end{eqnarray}
Introducing variables $L_i$ as before, we can rewrite the main term above  as
\begin{eqnarray}
    \label{eq:step3}
    \\
    \frac{(-1)^{\frac{k(k-1)}{2}+kN}}{k!}\prod_{j=1}^{2k}\frac{d}{d\alpha_j} e^{-\tfrac
    N2\sum_{j=1}^{2k}\alpha_j}
     \Delta^2\bigg(\frac{d}{dL}\bigg)
    \prod_{i=1}^k \bigg(\frac{d}{dL_i}\bigg)
    \prod_{i=1}^k  \bigg(\frac{1}{2\pi i}\int_{|w|=1}
    \frac{e^{L_i w }}{\prod_{j=1}^{2k}
    (w-\alpha_j)}\frac{dw}{w}\bigg). \notag
\end{eqnarray}
Now, by Lemma 6, the integral is
\begin{eqnarray}
    \int_
    {\sum_{j=1}^{2k}x_{ j}\le L_i}e^{\sum_{j=1}^{2k}x_{ j} \alpha_j}
    \prod_{  1\le j\le 2k}
    dx_{ j}.
\end{eqnarray}
Letting the variables in the $i$th integral
be $x_{i,j}$ we may express
the product of the $k$ integrals   as
\begin{eqnarray}
    \int_
    {\sum_{j=1}^{2k}x_{1,j}\le L_1} \dots
    \int_
    {\sum_{j=1}^{2k}x_{k,j}\le L_k}  e^{\sum_{i=1}^k\sum_{j=1}^{2k}x_{i,j}
    \alpha_j}
    \prod_{1\le i\le k\atop 1\le j\le 2k}
    dx_{i,j}.
\end{eqnarray}
We incorporate the factor $e^{-\tfrac N2\sum_{j=1}^{2k}\alpha_j}$ into
this product
and have
\begin{eqnarray}
    \int_
    {\sum_{j=1}^{2k}x_{1,j}\le L_1} \dots
    \int_
    {\sum_{j=1}^{2k}x_{k,j}\le L_k}
    e^{\sum_{j=1}^{2k} \alpha_j \big(\sum_{i=1}^{k}x_{i,j}-  N/ 2\big) }
    \prod_{1\le i\le k\atop 1\le j\le 2k}
    dx_{i,j}.
\end{eqnarray}
We differentiate this product of integrals with respect to each $\alpha_j$
and set each
$\alpha_j$ equal to 0 yielding
\begin{eqnarray}
    \int_
    {\sum_{j=1}^{2k}x_{1,j}\le L_1} \dots
    \int_
    {\sum_{j=1}^{2k}x_{k,j}\le L_k} \prod_{j=1}^{2k}\bigg(\sum_{i=1}^kx_{i,j}-
     \frac N 2\bigg)
    \prod_{1\le i\le k\atop 1\le j\le 2k}
    dx_{i,j}.
\end{eqnarray}
We want to compute this integral by multiplying out the product and using
Lemma 7. A good way to think about this is as follows.
By equation~(\ref{eq:multinom})
\begin{eqnarray}
    (A_1+\dots+A_k-A)^{2k}=\sum_{m\in P_O^{k+1}(2k)}
    \binom{2k}{m} (-A)^{m_0}A_1^{m_1}\cdots A_k^{m_k}.
\end{eqnarray}
When we multiply out the product we will have a
sum of $(k+1)^{2k}$ terms, each term being
a product of some number of factors $(-N/2)$
and $x_{i,j}$. Let $m\in P_O^{k+1}(2k)$
represent a generic term in which  $(-N/2)$ appears $m_0$ times,
and  factors $x_{1,j}$ appear for $m_1$ values of $j$,
and $x_{2,j}$ for $m_2$ values of $j$ and so on.
When we apply Lemma 7 to this term, when we perform the integration
over the variables $x_{1,j}$ the answer is solely determined by
$m_1$, the number of different $x_{1,j}$ that appear in this term.
Therefore, we find that
the product of integrals evaluates as
\begin{eqnarray}
    \sum_{m\in P_O^{k+1}(2k)} \binom{2k}{m}
     \left(-\frac{N}{2}\right)^{m_0}
    \frac{L_1^{2k+m_1}}{(2k+m_1)!}\cdots
    \frac{L_k^{2k+m_k}}{(2k+m_k)!}.
\end{eqnarray}

We now have that the quantity in equation~(\ref{eq:step3}) is equal to
\begin{eqnarray}
     \\
     \frac{(-1)^{\frac{k(k-1)}{2}+kN}}{k!}
     \Delta^2\bigg(\frac{d}{dL}\bigg)
     \prod_{i=1}^k \bigg(\frac{d}{dL_i}\bigg)
     \sum_{m\in P_O^{k+1}(2k)} \binom{2k}{m}
     \left(-\frac{N}{2}\right)^{m_0}
     \frac{L_1^{2k+m_1}}{(2k+m_1)!}\cdots
     \frac{L_k^{2k+m_k}}{(2k+m_k)!}. \notag
\end{eqnarray}
Now we need to carry out the differentiations with respect to the
$L_i$ and set the $L_i$ equal to $N$.
We perform the differentiations $\prod_{i=1}^k d/{dL_i}$ and obtain
\begin{eqnarray}
    \\
    \frac{(-1)^{\frac{k(k-1)}{2}+kN}}{k!}
    \Delta^2\bigg(\frac{d}{dL}\bigg)
    \sum_{m\in P_O^{k+1}(2k)} \binom{2k}{m}
    \left(-\frac{N}{2}\right)^{m_0}
    \frac{L_1^{2k-1+m_1}}{(2k-1+m_1)!}\cdots
    \frac{L_k^{2k+m_k}}{(2k+m_k)!}. \notag
\end{eqnarray}
Now the sum over $m_1,\dots m_k$ is a symmetric function of
the variables $L_i$.  Therefore, we can apply the second part
of Lemma~\ref{lemma:two vandermondes} to obtain that the above,
evaluated at $L_i=N$ is
\begin{eqnarray}
    &&\int_{U(N)} |\mathcal Z_A'(1)|^{2k} dA_N \\
    &&=
    (-1)^{\frac{k(k+1)}{2}}
    \sum_{m\in P_O^{k+1}(2k)} \binom{2k}{m}
    \left(-\frac{N}{2}\right)^{m_0}
    \det_{k\times k}\bigg( \frac{N^{2k+1+m_i-i-j}}{(2k+1+m_i-i-j)!}\bigg)+O(N^{k^2+2k-1})\notag
\end{eqnarray}
which we rewrite as
\begin{eqnarray}
    \\
    (-1)^{\frac{k(k+1)}{2}}
    N^{k^2+2k}\sum_{m\in P_O^{k+1}(2k)} \binom{2k}{m}
   (- 2)^{-m_0}
    \det_{k\times k}\bigg( \frac{1}{(2k+1+m_i-i-j)!}\bigg) +O(N^{k^2+2k-1}). \notag
\end{eqnarray}
Here the signs work out as in (\ref{eq: Z' as contour}). We factor
$1/(2k-i+m_i)!$ out of the $i$th row. The remaining determinant
has $i$th row
\begin{eqnarray}
    \\
    1, ~ (2k-i+m_i), ~ (2k-i+m_i)(2k-i-1+m_i), \quad \dots,
    \quad \prod_{j=1}^{k-1}(2k-i-j+1+m_i) \notag
\end{eqnarray}
This determinant is a polynomial in the $m_i$ of degree $0+1+\dots
+(k-1)=k(k-1)/2$ which vanishes whenever $m_j-m_i=j-i$; moreover
the part of it with degree $k(k-1)/2$ is precisely
$\Delta(m_1,\dots,m_k)=\prod_{1\le i< j\le k}(m_j-m_i)$.
Consequently the determinant evaluates to
\begin{eqnarray}
    \prod_{1\le i< j\le k}(m_j-m_i-j+i).
\end{eqnarray}
This concludes the evaluation of $b_k'$.

\vfill
\eject

\section{Numerical evaluation of $b_k$ and $b_k'$}

We have the following values for $b_k$:

\begin{eqnarray*}
b_1&=&\frac{1}{3}\end{eqnarray*} \begin{eqnarray*}
b_2&=& \frac{61}{2^5\cdot 3^2\cdot 5\cdot 7}\end{eqnarray*}
\begin{eqnarray*}
b_3&=& \frac{277}{2^{7}\cdot 3^4\cdot 5^2\cdot 7^2\cdot 11}\end{eqnarray*}
\begin{eqnarray*}
b_4&=& \frac{2275447}{2^{18}\cdot 3^{10}\cdot 5^4\cdot 7^3\cdot 11\cdot
13}\end{eqnarray*} \begin{eqnarray*}
b_5&=& \frac{3700752773}{2^{26}\cdot 3^{14}\cdot 5^6\cdot 7^4\cdot
11^2\cdot 13^2\cdot 17
\cdot 19}\end{eqnarray*} \begin{eqnarray*}
b_6&=& \frac{3654712923689}{2^{39}\cdot 3^{19}\cdot 5^9\cdot 7^6\cdot
11^3\cdot 13^3\cdot 17
\cdot 19\cdot 23}\end{eqnarray*} \begin{eqnarray*}
b_7&=& \frac{53\cdot 13008618017\cdot 143537}{2^{50}\cdot 3^{28}\cdot
5^{13}\cdot 7^8\cdot 11^5\cdot 13^4\cdot 17^2
\cdot 19^2\cdot 23}\end{eqnarray*} \begin{eqnarray*}
b_8&=& \frac{41\cdot 359\cdot 5505609492791\cdot 3637}{2^{68}\cdot
3^{35}\cdot 5^{16}\cdot 7^{11}\cdot 11^6\cdot 13^5\cdot 17^3
\cdot 19^2\cdot 23\cdot 29 \cdot 31}\end{eqnarray*} \begin{eqnarray*}
b_9&=& \frac{757\cdot 45742439\cdot 60588179\cdot 13723}{2^{84}\cdot
3^{42}\cdot 5^{21}\cdot 7^{14}\cdot 11^8\cdot 13^6\cdot 17^4
\cdot 19^3\cdot 23^2\cdot 29 \cdot 31}\end{eqnarray*} \begin{eqnarray*}
b_{10}&=& \frac{652071900673\cdot 241845775551409}{2^{105}\cdot
3^{55}\cdot 5^{25}\cdot 7^{17}\cdot 11^{10}\cdot 13^8\cdot 17^5
\cdot 19^4\cdot 23^3\cdot 29  \cdot 37}\end{eqnarray*} \begin{eqnarray*}
b_{11}&=& \frac{1318985497\cdot 578601141598041214011811}{2^{121}\cdot
3^{64}\cdot 5^{31}\cdot 7^{19}\cdot 11^{12}\cdot 13^{9}\cdot 17^7
\cdot 19^6\cdot 23^4\cdot 29^2 \cdot 31^2\cdot 37\cdot 41\cdot
43}\end{eqnarray*} \begin{eqnarray*}
b_{12}&=& \frac{113\cdot 206489633386447920175141\cdot 51839\cdot 14831}
{2^{150}\cdot 3^{75}\cdot 5^{37}\cdot 7^{23}\cdot 11^{15}\cdot
13^{12}\cdot 17^7
\cdot 19^7\cdot 23^5\cdot 29^3 \cdot 31^2 \cdot 37 \cdot 41\cdot 43\cdot
47}\end{eqnarray*} \begin{eqnarray*}
b_{13}&=& \frac{4670754069404622871904068067089635254838677}{2^{174}\cdot
3^{90}\cdot 5^{42}\cdot 7^{28}\cdot 11^{17}
\cdot 13^{14}\cdot 17^{10}
\cdot 19^9\cdot 23^6\cdot 29^3 \cdot 31^3\cdot 37^2 \cdot 41\cdot 43\cdot
47}\end{eqnarray*} \begin{eqnarray*}
b_{14}&=& \frac{107\cdot 194946046688455595346779341\cdot
996075171809335069}{2^{203}\cdot 3^{103}\cdot 5^{50}\cdot 7^{31}\cdot
11^{20}
\cdot 13^{17}\cdot 17^{12}
\cdot 19^{10}\cdot 23^7\cdot 29^4 \cdot 31^4\cdot 37^2 \cdot 41\cdot
43\cdot 47\cdot 53}\end{eqnarray*} \begin{eqnarray*}
b_{15}&=& \frac{29547975377\cdot 3981541\cdot
1807995588661527603489333681461\cdot 1584311}{2^{230}\cdot 3^{117}\cdot
5^{57}
\cdot 7^{37}\cdot 11^{22}
\cdot 13^{19}\cdot 17^{14}
\cdot 19^{12}\cdot 23^9\cdot 29^5 \cdot 31^5\cdot 37^3
\cdot 41^2\cdot 43^2\cdot 47\cdot 53\cdot 59}\end{eqnarray*}

\vfill
\eject

We have the following values for $b_k'$:

\begin{eqnarray*}
b_1'&=&\frac{1}{2^2\cdot 3}\end{eqnarray*} \begin{eqnarray*}
b_2'&=& \frac{1}{2^6\cdot 3\cdot 5\cdot 7}\end{eqnarray*}
\begin{eqnarray*}
b_3'&=& \frac{1}{2^{12}\cdot 3^2\cdot 5^2\cdot 7^2\cdot 11}\end{eqnarray*}
\begin{eqnarray*}
b_4'&=& \frac{31}{2^{20}\cdot 3^{10}\cdot 5^4\cdot 7^2\cdot 11\cdot
13}\end{eqnarray*} \begin{eqnarray*}
b_5'&=& \frac{227}{2^{30}\cdot 3^{12}\cdot 5^6\cdot 7^4\cdot 11\cdot
13^2\cdot 17
\cdot 19}\end{eqnarray*} \begin{eqnarray*}
b_6'&=& \frac{67 \cdot 1999}{2^{42}\cdot 3^{19}\cdot 5^9\cdot 7^6\cdot
11^3\cdot 13^3\cdot 17
\cdot 19\cdot 23}\end{eqnarray*} \begin{eqnarray*}
b_7'&=& \frac{43\cdot  46663}{2^{56}\cdot 3^{28}\cdot 5^{13}\cdot 7^8\cdot
11^4\cdot 13^3\cdot 17^2
\cdot 19^2\cdot 23}\end{eqnarray*} \begin{eqnarray*}
b_8'&=& \frac{46743947}{2^{72}\cdot 3^{34}\cdot 5^{16}\cdot 7^{11}\cdot
11^6\cdot 13^4\cdot 17^3
\cdot 19^2\cdot 23\cdot 29 \cdot 31}\end{eqnarray*} \begin{eqnarray*}
b_9'&=& \frac{19583\cdot16249}{2^{90}\cdot 3^{42}\cdot 5^{21}\cdot
7^{14}\cdot 11^8\cdot 13^6\cdot 17^3
\cdot 19^3\cdot 23^2\cdot 29 \cdot 31}\end{eqnarray*} \begin{eqnarray*}
b_{10}'&=& \frac{3156627824489}{2^{110}\cdot 3^{55}\cdot 5^{25}\cdot
7^{17}\cdot 11^{10}\cdot 13^8\cdot 17^5
\cdot 19^4\cdot 23^3\cdot 29 \cdot 31\cdot 37}\end{eqnarray*}
\begin{eqnarray*}
b_{11}'&=& \frac{59\cdot 11332613\cdot 33391}{2^{132}\cdot 3^{63}\cdot
5^{31}\cdot 7^{18}\cdot 11^{12}\cdot 13^{10}\cdot 17^5
\cdot 19^5\cdot 23^4\cdot 29^2 \cdot 31^2\cdot 37\cdot 41\cdot
43}\end{eqnarray*} \begin{eqnarray*}
b_{12}'&=& \frac{241\cdot 251799899121593}{2^{156}\cdot 3^{75}\cdot
5^{37}\cdot 7^{23}\cdot 11^{15}\cdot 13^{12}\cdot 17^8
\cdot 19^7\cdot 23^4\cdot 29^3 \cdot 31^2 \cdot 41\cdot 43\cdot
47}\end{eqnarray*} \begin{eqnarray*}
b_{13}'&=& \frac{285533 \cdot 37408704134429 }{2^{182}\cdot 3^{90}\cdot
5^{42}\cdot 7^{28}\cdot 11^{17}\cdot 13^{14}\cdot 17^{10}
\cdot 19^8\cdot 23^5\cdot 29^3 \cdot 31^3\cdot 37^2
\cdot 41\cdot 43\cdot 47}\end{eqnarray*}
\begin{eqnarray*}
b_{14}' &=& \frac{197\cdot 1462253323\cdot 6616773091 }{2^{210}
\cdot 3^{100}\cdot 5^{50}\cdot 7^{31}\cdot 11^{20}\cdot 13^{17}\cdot
17^{12}
\cdot 19^{10}\cdot 23^7\cdot 29^4 \cdot 31^4\cdot 37^2
\cdot 41\cdot 43\cdot 47\cdot 53}\end{eqnarray*}
\begin{eqnarray*}
b_{15}' &=& \frac{1625537582517468726519545837}{2^{240}
\cdot 3^{117}\cdot 5^{57}\cdot 7^{37}\cdot 11^{22}\cdot 13^{19}\cdot
17^{14}
\cdot 19^{11}\cdot 23^9\cdot 29^5 \cdot 31^5\cdot 37^3
\cdot 41^2\cdot 43^2\cdot 47\cdot 53\cdot 59}\end{eqnarray*}

 \vfill
 \eject

\section{Acknowledgements}

The authors are grateful to AIM and the Isaac Newton Institute for
very generous support and hospitality. JBC was supported by
the NSF, MOR by the NSF and NSERC, and NCS by an EPSRC Advanced Research 
Fellowship.

\end{document}